\documentclass[11pt]{article}

\usepackage{graphicx,enumitem}
\usepackage{graphicx,psfrag,epsf}
{%
\centering\addtocounter{figure}{1}% if caption at bottom
\begin{enumerate}[%
itemsep=2pt,parsep=0em,
label={(\alph*)},
ref={\thefigure.(\alph*)}
]}%
{\end{enumerate}\addtocounter{figure}{-1}}

\usepackage[margin=1in]{geometry}

\usepackage[utf8]{inputenc}
\usepackage{amsmath}
\usepackage[english]{babel}
\usepackage{amssymb}
\usepackage{amsfonts}
\usepackage{bm}
\usepackage{graphicx}
\usepackage{amsmath}
\usepackage{amsthm}
\usepackage{xcolor}
\usepackage[hidelinks=true]{hyperref}
\usepackage{parskip}
\usepackage{multicol}
\usepackage{float}

\usepackage{verbatim}

\floatstyle{plaintop}
\restylefloat{table}
\usepackage[tableposition=top]{caption}

\usepackage{rotating}
\usepackage{subcaption}
\usepackage{multirow}

\usepackage{chngcntr}
\usepackage{apptools}
\AtAppendix{\counterwithin{lem}{subsection}}
\AtAppendix{\counterwithin{proposition}{subsection}}
\AtAppendix{\counterwithin{equation}{subsection}}

 \theoremstyle{definition}
\newtheorem{thm}{Theorem}[section]

\newtheorem{lem}{Lemma}[section]

\newtheorem{mydef}{Definition}[section]

\usepackage{booktabs}
\usepackage{animate}

\usepackage{natbib}
\hypersetup{urlcolor=blue, citecolor=blue, colorlinks=true, linkcolor=blue} 

\numberwithin{equation}{section}
\numberwithin{figure}{section}
\numberwithin{table}{section}

\newcommand{\blind}{1}

\usepackage{blindtext}
\usepackage{lineno}

\begin{document}

\def\spacingset#1{\renewcommand{\baselinestretch}%
{#1}\small\normalsize} \spacingset{1}

\if1\blind
{
 \title{\textbf{Covariance Functions for Multivariate Gaussian Fields Evolving Temporally Over Planet Earth}}

\author{Alfredo Alegr\'ia$^{1}$, Emilio Porcu$^{1,2}$, Reinhard Furrer$^{3}$ and Jorge Mateu$^{4}$\\\\
 $^{1}$School of Mathematics \& Statistics, Newcastle University,\\  Newcastle Upon Tyne, United Kingdom.\\ 
 $^{2}$Departamento de Matem\'atica, Universidad T\'ecnica Federico Santa Maria,\\  Valpara\'iso, Chile.\\ 
 $^{3}$Department of Mathematics and Department of Computer Science,\\ University of Zurich, Zurich, Switzerland.\\
$^{4}$Department of Mathematics, University Jaume I,\\ Castell\'on, Spain.}

 \maketitle
} \fi

\if0\blind
{
  \bigskip

 % \begin{center}
  \title{\textbf{Covariance Functions for Multivariate Gaussian Fields evolving temporally over Planet Earth}}
   \maketitle
%\end{center}
  \medskip
} \fi

\begin{abstract}
\noindent The construction of valid and flexible cross-covariance functions is a fundamental task for modeling  multivariate  space-time data arising from  climatological and oceanographical   phenomena. Indeed, a suitable specification of the covariance structure allows to capture both the  space-time  dependencies between the observations and  the development of accurate predictions. For data observed over large portions of planet Earth it is necessary to take into account the curvature of the planet. Hence the need for random field models defined over spheres across time.  In particular, the associated covariance function should depend on the geodesic distance, which is the most natural metric over the spherical surface.      
 In this work,  we  propose a flexible parametric family of matrix-valued covariance functions, with both marginal and cross structure being of  the   Gneiting type.  We additionally introduce a different multivariate Gneiting model   based on the adaptation of the latent dimension approach to the spherical context.   Finally, we assess the performance of our models through the study of  a  bivariate space-time data set of surface air temperatures and precipitations.  \\

\end{abstract}

\spacingset{1} % DON'T change the spacing!

\def\bY{\bm{\Upsilon}}
\def\bC{\bm{C}}
\def\bZ{\bm{Z}}
\def\R{\mathbb{R}}

\section{Introduction}

Monitoring several georeferenced variables is a common practice in a wide range of disciplines such as climatology and oceanography. The phenomena under study are often  observed over  large portions of the Earth and     across several time instants.  
Since  there is only a finite sample from the involved variables, geostatistical models  are a useful tool to capture both the  spatial and temporal  interactions between the observed   data, as well as the uncertainty associated to the limited available information  (\citealp{Cressie:1993}; \citealp{Wackernagel:2003};  \citealp{Gneiting:Genton:Guttorp:2007}).    
The geostatistical approach consists in modeling the observations  as a partial realization of a space-time multivariate random field (RF), denoted as  $\{\bm{Z}(\bm{x},t)=(Z_1(\bm{x},t),\hdots,Z_m(\bm{x},t))^\top: (\bm{x},t)\in \mathcal{D}\times \mathcal{T} \}$,   where $\top$ is the transpose operator and $m$ is  a positive integer representing the number of components of the field.   If $m=1$, we say that $Z(\bm{x},t)$ is a univariate (or scalar-valued) RF, whereas   for $m > 1$,  $\bm{Z}(\bm{x},t)$ is called an $m$-variate (or vector-valued) RF. Here, $\mathcal{D}$ and $\mathcal{T}$ denote the spatial and temporal domains, respectively.    Throughout, we assume that $\bm{Z}(\bm{x},t)$ is a zero mean Gaussian field,  so that  a suitable specification of its covariance structure  is crucial to develop both accurate inferences and predictions over unobserved sites  \citep{Cressie:1993}.

 Parametric families of matrix-valued covariance functions are typically given in terms of  Euclidean distances.   The literature for this case is extensive and we refer the reader to the review  by \cite{Genton:Kleiber:2014} with the references therein. The main motivation to consider the Euclidean metric is the existence of several methods for projecting the geographical coordinates, longitude and latitude, onto the plane. 
  However, when a phenomenon  is observed over large portions of  planet Earth,  the approach based on projections  generates   distortions in the distances associated to distant  locations on the globe. The reader is referred to  \cite{BIOM:BIOM040302}, where  the impact of the different types  of  projections with respect to spatial inference is discussed. 
  
   Indeed,  the geometry of the Earth must be considered. Thus, it is more realistic  to work under the framework of RFs defined spatially on a sphere (see \citealp{marinucci2011random}).   Let $d$ be a positive integer. The $d$-dimensional unit sphere is denoted as $\mathbb{S}^d:=\{\bm{x}\in\mathbb{R}^{d+1}: \|\bm{x}\| = 1\}$, where $\|\cdot\|$ represents the Euclidean norm.
The most accurate metric in the spherical scenario is the  geodesic (or great circle) distance, which     roughly speaking corresponds   to the arc joining any two  points located on the sphere, measured along a path on the spherical surface. Formally, the geodesic distance is defined as 
the mapping $\theta:\mathbb{S}^d \times \mathbb{S}^d \rightarrow [0,\pi]$ given by $\theta := \theta(\bm{x},\bm{y}) = \arccos   (\bm{x}^\top\bm{y})$.

The construction of valid and flexible parametric covariance functions in terms of the geodesic distance is a challenging problem and requires the application of the theory of positive definite functions on spheres (\citealp{schoenberg1942};  \citealp{Yaglom:1987};  \citealp{hannan2009multiple}; \citealp{berg2016schoenberg}).    
In the univariate   and merely spatial case, \cite{Huang:Zhang:Robenson:2011} study the validity of some specific  covariance functions. The essay by \cite{gneiting2013} contains a wealth of  results related to the validity of a wide range of covariance families.  Other related works are the study of star-shaped random particles \citep{hansen2011levy} and convolution roots \citep{ziegel2014convolution}.
However,   the   spatial and spatio-temporal covariances  in the multivariate case are  still unexplored, with the work of  \cite{PBG16} being a notable exception.

The Gneiting class \citep{doi:10.1198/016214502760047113}  is one of the most popular space-time covariance families and some adaptations in terms of geodesic distance have been given by \cite{PBG16}.  In this paper, we extend to the multivariate scenario the modified Gneiting class introduced by \cite{PBG16}.  Furthermore, we adapt the latent dimension approach   to the spherical context  \citep{Apanasovich:Genton:2010} and  we then generate  additional Gneiting type  matrix-valued covariances.   The proposed models are non-separable with respect to the components of the field nor with respect to the space-time interactions.  To obtain these results, we have demonstrated several technical results that can be useful to develop new research in this area.  Our findings are illustrated through a    real data application.  In particular, we analyze a bivariate space-time data set of surface air temperatures and precipitations. These  variables have been obtained from the  Community Climate System Model (CCSM) provided by NCAR (National Center for Atmospheric Research).

 The remainder of the article is organized as follows.  Section \ref{preliminares} introduces some concepts and construction principles of cross-covariances on spheres across time. In Section \ref{gneiting_families} we propose some multivariate Gneiting type covariance families.  Section \ref{sec_data}  contains    a real data example of surface air temperatures and precipitations. Section \ref{discussion} concludes the paper with a discussion.

\section{Fundaments and Principles of   Space-time Matrix-valued Covariances}
\label{preliminares}

This section is devoted to introduce the basic material related to multivariate fields on spheres across time.  Let $\bm{Z}(\bm{x},t)$, for $\bm{x}\in\mathbb{S}^d$ and $t\in\R$, be an $m$-variate space-time field and let  $\textbf{\text{C}}: [0,\pi]\times \mathbb{R}\rightarrow \mathbb{R}^{m\times m}$ be a continuous matrix-valued mapping, whose elements are defined as $\text{C}_{ij}(\theta,u) = \text{cov}\{Z_i(\bm{x},t+u),Z_j(\bm{y},t)\}$, where $\theta=\theta(\bm{x},\bm{y})$. Following \cite{PBG16}, we say that $\textbf{C}$ is geodesically isotropic in space and  stationary in time.  Throughout, the diagonal elements $\text{C}_{ii}$ are called marginal covariances, whereas the off-diagonal members $\text{C}_{ij}$ are called cross-covariances. Any  parametric representation of  $\textbf{\text{C}}$ must respect the positive definite condition, i.e. it must satisfy 
\begin{equation} \label{posdef}
   \sum_{\ell=1}^{n} \sum_{r=1}^n   \bm{a}_{\ell}^\top \textbf{\text{C}}(\theta(\bm{x}_\ell,\bm{x}_r),t_\ell-t_r) \bm{a}_{r} \geq 0,
 \end{equation}
for all positive integer $n$, $\{(\bm{x}_1,t_1),\hdots,(\bm{x}_n,t_n)\}\subset \mathbb{S}^d\times\mathbb{R}$ and $\{\bm{a}_1,\hdots,\bm{a}_n\}\subset \mathbb{R}^m$.  Appendix \ref{background} contains a complete background material on matrix-valued positive definite functions on Euclidean and spherical domains.

We call the mapping $\textbf{\text{C}}$  {\em space-time $m$-separable} if there exists two mappings $\textbf{\text{C}}_{{\tt S}}:[0,\pi]\rightarrow \mathbb{R}^{m\times m}$ and $\textbf{\text{C}}_{{\tt T}}:\mathbb{R}\rightarrow \mathbb{R}^{m\times m}$, being  merely spatial and temporal matrix-valued covariances, respectively, such that
$$ \textbf{\text{C}}(\theta,u)= \textbf{\text{C}}_{{\tt S}}(\theta) \circ \textbf{\text{C}}_{{\tt T}}(u), \qquad (\theta,u) \in [0,\pi] \times \R, $$
where $\circ$ denotes the Hadamard product.   We call the space-time $m$-separability property {\em complete} if there exists a symmetric, positive definite matrix $\textbf{\text{A}} \in \mathbb{R}^{m \times m}$, a univariate spatial covariance $\text{C}_{{\tt S}}: [0,\pi]\rightarrow \mathbb{R}$ , and a univariate temporal covariance $\text{C}_{{\tt T}}: \mathbb{R}\rightarrow \mathbb{R}$, such that 
$$ \textbf{\text{C}}(\theta,u) = \textbf{\text{A}} \text{C}_{{\tt S}}(\theta) \text{C}_{{\tt T}}(u), \qquad (\theta,u) \in [0,\pi] \times \R. $$  
We finally call {\em $m$-separable} the mapping $\textbf{\text{C}}$ for which there exists a univariate space-time covariance $\text{C}:[0,\pi]\times \mathbb{R}\rightarrow \mathbb{R}$, and a matrix $\textbf{\text{A}}$, as previously defined, such that 
$$ \textbf{\text{C}}(\theta,u) = \textbf{\text{A}} \text{C}(\theta,u), \qquad (\theta,u) \in [0,\pi] \times \R, $$
and of course the special case $\text{C}(\theta,u)= \text{C}_{{\tt S}}(\theta) \text{C}_{{\tt T}}(u)$ offers complete space-time $m$-separability as previously discussed. 

Separability is a very useful property for both modeling and estimation purposes, because the related covariance matrices admit nice factorizations, with consequent alleviation of the computational burdens. At the same time, it is a very unrealistic assumption and the literature has been focussed on how to develop non-separable models.   How to escape from separability is a major deal, and we list some strategies that can be adapted from others proposed in Euclidean spaces. To the knowledge of the authors, none of these strategies have ever been implemented on spheres or spheres across time.

\begin{description}

\item[Linear Models of Coregionalization]
Let $q$ be a positive integer. Given a collection of matrices $\textbf{\text{A}}_{k}$, $k=1,\ldots,q$, and a collection of univariate space-time covariances $\text{C}_k :[0,\pi]\times \mathbb{R}\rightarrow \mathbb{R}$, the linear model of coregionalization (LMC) has expression
\begin{equation*}\label{LMC}
 \textbf{\text{C}}(\theta,u)= \sum_{k=1}^q \textbf{\text{A}}_k  \text{C}_k(\theta,u), \qquad (\theta,u) \in [0,\pi] \times \R, 
 \end{equation*}
where a simplification of the type $\text{C}_k(\theta,u)= \text{C}_{k,{\tt S}}(\theta) \text{C}_{{\tt T}}(u)$ might be imposed. This model has several drawbacks that have been discussed in \cite{Gneiting:Kleibler:Schlather:2010} as well as in \cite{daley2015}.

\item[Lagrangian Frameworks] Let $\boldsymbol{Z}$ be an $m$-variate Gaussian field on the sphere with covariance $\textbf{\text{C}}_{{\tt S}}: [0,\pi]\rightarrow \mathbb{R}^{m\times m}$. Let ${\cal R}$ be a random orthogonal matrix with a given  probability law. Let 
$$ \boldsymbol{Y}(\boldsymbol{x},t):=   \boldsymbol{Z} \left ( {\cal R}^{t} \boldsymbol{x} \right ), \qquad (\boldsymbol{x},t) \in \mathbb{S}^{d} \times \R. $$
Then, $\boldsymbol{Y}$ is a RF with transport effect over the sphere. Clearly, $\boldsymbol{Y}$ is not Gaussian and evaluation of the corresponding covariance might be cumbersome, as shown in \cite{alegria2016dimple}.

\item[Multivariate Parametric Adaptation] Let $p$ be a positive integer. Let $\text{C}(\cdot,\cdot; \bm{\lambda})$, for $\bm{\lambda} \in \R^p$, be a univariate space-time covariance. Let $\boldsymbol{\lambda}_{ij} \in \R^p$, for $i,j=1,\ldots,m$. For $|\rho_{ij}|\le 1$ and  $\rho_{ii}=1$,  and $\sigma_{ii}>0$, find the parametric restrictions such that $\textbf{\text{C}}:[0,\pi] \times \R \rightarrow \R^{m\times m}$ defined through
$$ \text{C}_{ij}(\theta,u) = \sigma_{ii}\sigma_{jj}\rho_{ij} \text{C}(\theta,u; \boldsymbol{\lambda}_{ij}), \qquad (\theta,u) \in [0,\pi] \times \R, $$ 
is a valid covariance. Sometimes the restriction on the parameters can be very severe, in particular when $m$ is bigger than 2. In Euclidean spaces this strategy has been adopted by \cite{Gneiting:Kleibler:Schlather:2010} and \cite{doi:10.1080/01621459.2011.643197} for the Mat{\'e}rn model, and by \cite{daley2015} for models with compact support. 

\item[Scale Mixtures] Let $(\Omega, {\cal A}, \mu) $ be a measure space. Let $\textbf{\text{C}}: [0,\pi ]\times \R \times \Omega \rightarrow \R^{m\times m}$ such that 
\begin{enumerate}
\item $\textbf{\text{C}}(\cdot,\cdot, \xi)$ is a valid covariance for all $\xi$ in $\Omega$;
\item $\textbf{\text{C}}(\theta,u, \cdot) \in L_1 \left ( \Omega, {\cal A}, \mu \right ) $ for all $(\theta,u) \in [0,\pi] \times \R$.
\end{enumerate}
Then,
$$ \int_{\Omega} \textbf{\text{C}}(\theta,u,\xi) \mu({\rm d} \xi) $$
is still a valid covariance \citep{Porcu20111293}. Of course, simple strategies can be very effective. For instance, one might assume that $\textbf{\text{C}}(\cdot,\cdot,\xi)= \text{C}(\cdot,\cdot) \textbf{\text{A}}(\xi)$, with $\text{C}$ a univariate covariance and $\textbf{\text{A}}(\xi) \in \R^{m \times m}$ being a positive definite matrix for any $\xi \in \Omega$ and such that the hypothesis of integrability above is satisfied. 
 \end{description}

\section{Matrix-valued Covariance Functions of Gneiting Type}
\label{gneiting_families}

This section provides some general results for the construction of  multivariate   covariances  for fields defined on $\mathbb{S}^d\times \R$. The most important feature of the proposed models   is that they are  non-separable mappings,  allowing more flexibility in the study of space-time phenomena. In particular, we focus on covariances with  a Gneiting's structure \citep{doi:10.1198/016214502760047113}.   See also \cite{bourotte2016flexible} for multivariate  Gneiting classes on Euclidean spaces. 

Recently,  \cite{PBG16} have proposed some versions of the Gneiting model  for  RFs with spatial locations on the unit sphere. Let us provide a brief review.     The mapping $\text{C}:[0,\pi]\times \mathbb{R}\rightarrow \mathbb{R}$ defined as
\begin{equation}
\label{gneiting_porcu}
\text{C}(\theta,u) =  \frac{1}{    \left\{  f(\theta) |_{[0,\pi]}  \right\}^{1/2}      } g\left(   \frac{|u|^2}{   f(\theta) |_{[0,\pi]}    } \right), \qquad \theta \in [0,\pi], u \in \mathbb{R},
\end{equation}
is called  Gneiting model. Arguments in \cite{PBG16} show that $\text{C}$  is a valid univariate covariance,  for all $d\in\mathbb{N}$, 
for $g:[0,\infty)\rightarrow [0,\infty)$ being a completely monotone function, i.e.  $g$ is infinitely differentiable on $(0,\infty)$ and $(-1)^n g^{(n)}(t)\geq 0$, for all $n\in\mathbb{N}$ and $t> 0$, and $f:[0,\infty)\rightarrow (0,\infty)$ is  a positive function with completely monotone derivative. Such functions $f$ are called Bernstein functions  \citep{porcu2011}.  Here, $f(\theta)|_{[0,\pi]}$ denotes the restriction of the mapping $f$ to the interval $[0,\pi]$.   Tables \ref{cm} and \ref{bernstein}  contain some examples of completely monotone and Bernstein functions, respectively. Additional properties about these functions are studied in \cite{porcu2011}. 
 
We also pay attention to the modified Gneiting class \citep{PBG16} defined through 
 \begin{equation}
 \label{mod_uni}
 \text{C}(\theta,u) = \frac{ 1}{f(|u|)^{n+2}} g \left (\theta f(|u|) \right ), \qquad (\theta,u) \in  [0,\pi] \times \mathbb{R},
 \end{equation}
where $n\leq 3$ is a positive integer,  $g: [0,\infty)\rightarrow [0,\infty)$ is a completely monotone function and $f: [0,\infty)\rightarrow (0,\infty)$ is a strictly increasing and concave function on the positive real line.  The mapping (\ref{mod_uni}) is a valid covariance for any $d\leq 2n+1$.

\begin{table}
\centering{
\caption{Some examples of completely monotone functions. Here, $K_\nu$ denotes the modified Bessel function of second kind of degree $\nu$.}
\label{cm}
\begin{tabular}{cc}\hline \hline
Function &  Parameters restriction\\ \hline
$g(t) = \exp(-ct^\gamma)$ & $c>0$, $0< \gamma \leq 1$\\
$g(t) = (2^{\nu-1}\Gamma(\nu))^{-1}(c\sqrt{t})^\nu K_\nu(c\sqrt{t})$ & $c>0$, $\nu>0$\\
$g(t) = (1+ct^\gamma)^{-\nu}$ & $c>0$, $0< \gamma \leq 1$, $\nu>0$\\
$g(t) = 2^\nu (\exp(c\sqrt{t})+\exp(-c\sqrt{t}))^{-\nu}$ & $c>0$, $\nu>0$\\ \hline 
\end{tabular}
}
\end{table}

\begin{table}
\centering{
\caption{Some examples of Bernstein  functions.}
\label{bernstein}
\begin{tabular}{cc}\hline \hline
Function &  Parameters restriction\\ \hline
$f(t) = (at^\alpha+1)^\beta$ &    $a>0$, $0<\alpha\leq 1$, $0\leq \beta\leq 1$ \\
$f(t) = \ln(at^\alpha+b)/\ln(b)$ &  $a>0$, $b>1$ $0<\alpha\leq 1$  \\
$f(t) = (at^\alpha+b)/(b(at^\alpha+1))$ &  $a>0$, $0< b \leq 1$, $0< \alpha \leq 1$  \\ \hline 
\end{tabular}
}
\end{table}

\subsection{Multivariate  Modified Gneiting Class on  Spheres Acros Time}

We are now able to illustrate the main result within this section. Specifically, we propose a multivariate space-time class generating Gneiting functions with different  scale parameters. 

\begin{thm} \label{los_tanVores}
Let $m \ge 2$ and $n \leq 3$ be positive integers. Let $g:[0,\infty) \to [0,\infty)$ be a completely monotone function. Consider $f: [0,\infty) \to (0,\infty)$ being strictly increasing and concave. Let $\sigma_{ii} > 0$, $|\rho_{ij}| \leq 1$ and $c_{ij}>0$, for $i,j=1,\ldots,m$, be constants yielding the  additional condition
\begin{equation}
\label{askey3}
\sum_{i\neq j} |\rho_{ij}|  (c_{ii}/c_{ij})^{n+1} \leq 1,
\end{equation}
for each $i=1,\hdots,m$. Then, the mapping $\textbf{\text{C}} $, whose members $\text{C}_{ij}:[0,\pi] \times \mathbb{R} \to \R$  are defined through
\begin{equation} \label{los_tanVores2}
\text{C}_{ij}(\theta,u) = \frac{  \sigma_{ii} \sigma_{jj} \rho_{ij}}{f(|u|)^{n+2}} g \left ( \frac{\theta f(|u|)}{c_{ij}} \right ), \qquad (\theta,u) \in  [0,\pi] \times \mathbb{R},
\end{equation}
is a matrix-valued covariance for any $d\leq 2n+1$. 
\end{thm}
The condition (\ref{askey3}) comes from the arguments in  \cite{daley2015}. The proof of Theorem \ref{los_tanVores} is deferred to Appendix \ref{modified_proof}, coupled with the preliminary notation and results introduced in Appendix \ref{background}.

For example, taking $n=1$,  $g(t) = \exp(-3t)$   and  $f(t) = 1+ 1.7  t/c_T$, for $c_T>0$, we can use  Equation (\ref{los_tanVores2}), and the restrictions of Theorem \ref{los_tanVores}, to generate a model of the form
\begin{equation}\label{modified1}
\text{C}_{ij}(\theta,u) = \frac{\sigma_{ii}\sigma_{jj} \rho_{ij}}{\left(1+ \frac{1.7   |u|}{c_T}\right)^{3}} \exp\left\{ - \frac{3\theta\left(1+ \frac{1.7  |u|}{c_T}\right)}{c_{ij}}\right\}, \qquad (\theta,u) \in  [0,\pi] \times \mathbb{R}.
\end{equation}
The special parameterization used in the covariance (\ref{modified1}) ensures that $\text{C}_{ij}(\theta,0)/(\sigma_{ii} \sigma_{jj}) < 0.05$ and $\text{C}_{ij}(0,|u|)/(\sigma_{ii} \sigma_{jj}) < 0.05$,  for $\theta>c_{ij}$ and $|u|>c_T$, respectively.

\subsection{A Multivariate Gneiting Family Based on  Latent Dimension Approaches}

We now  propose  families of  matrix-valued covariances  obtained from  univariate models.     The method used is known as  \textit{latent dimension} approach and has been studied in the Euclidean case by \cite{Porcu:Gregori:Mateu:2006},  \cite{Apanasovich:Genton:2010} and \cite{Porcu20111293}.   Next, we illustrate the spirit of this approach.

Consider  a univariate  Gaussian RF defined on the product space $\mathbb{S}^d \times \mathbb{R} \times \mathbb{R}^k$, for some positive integers $d$ and $k$, namely  $\{Y(\bm{x},t,\bm{\xi}): (\bm{x},t,\bm{\xi}) \in \mathbb{S}^d \times \mathbb{R} \times \mathbb{R}^k\}$.   Suppose that there exists a mapping $\text{K}:[0,\pi]\times \mathbb{R} \times  \mathbb{R}^k \rightarrow \mathbb{R}$ such that $\text{cov}\{Y(\bm{x},t,\bm{\xi}_1),Y(\bm{y},s,\bm{\xi}_2)\} = \text{K}(\theta(\bm{x},\bm{y}),t-s,\bm{\xi}_1-\bm{\xi}_2)$, for all $\bm{x},\bm{y}\in\mathbb{S}^d$, $t,s\in\mathbb{R}$ and $\bm{\xi}_1,\bm{\xi}_2\in\mathbb{R}^k$. The idea is to define the components of an $m$-variate RF  on $\mathbb{S}^d \times \mathbb{R}$ as
$$Z_i(\bm{x},t) = Y(\bm{x},t,\bm{\xi}_i),\qquad  \bm{x}\in\mathbb{S}^d, t\in\mathbb{R}, \bm{\xi}_i \in\mathbb{R}^k,$$
 for $i=1,\hdots,m$.   Thus, the resulting covariance, $\textbf{\text{C}}(\cdot)=[\text{C}_{ij}(\cdot)]_{i,j=1}^m$, associated to $\bm{Z}(\bm{x},t)$ is given by  $$\text{C}_{ij}(\theta,u) = \text{K}(\theta, u,\bm{\xi}_i-\bm{\xi}_j),  \qquad \theta\in[0,\pi], u\in\mathbb{R}.$$ 
Here, the vectors $\bm{\xi}_i$ are handled as additional parameters of the model. The following  theorem allows to construct different versions of the Gneiting model, using the latent dimension technique.

\begin{thm}
\label{prop_latent2}
Consider the  positive integers $d$, $k$ and $l$. Let $g$ be a completely monotone function and $f_i$, $i=1,2$,    Bernstein functions. Then,
\begin{equation}
\label{eq_latent2}
\text{K}(\theta,\bm{u},\bm{v}) =  \frac{1}{        \{f_2( \theta ) |_{[0,\pi]}\}^{l/2}    \left\{f_1\left[  \frac{\|\bm{v}\|^2}{  f_2(\theta) |_{[0,\pi]}} \right]  \right\}^{k/2}      } g\left(   \frac{\|\bm{u}\|^2}{   f _1\left[  \frac{\|\bm{v}\|^2}{f_2(\theta) |_{[0,\pi]}}  \right]    } \right), \qquad  (\theta,\bm{u}, \bm{v}) \in [0,\pi]\times \mathbb{R}^l\times \mathbb{R}^k,
\end{equation}
is a  univariate covariance for  fields defined  on $\mathbb{S}^d\times\R^l \times \R^k$.
\end{thm}
The proof of  Theorem \ref{prop_latent2} requires technical lemmas and is deferred to  Appendix  \ref{latent2}, coupled with the preliminary results introduced in Appendix \ref{background}. 

In order to avoid an excessive number of parameters, we follow the  parsimonious strategy proposed by \cite{Apanasovich:Genton:2010}.   Consider $k=1$ and the scalars $\{\xi_{1},\hdots,\xi_{m}\}$. We can consider the parameterization $\zeta_{ij}=|\xi_i-\xi_j|^2$, with $\zeta_{ij}=\zeta_{ji}>0$ and $\zeta_{ii}=0$, for all $i$ and $j$.   

Following  \cite{Apanasovich:Genton:2010},   a LMC based on the latent dimension approach can be used to achieve different marginal structures.  Indeed, suppose that the  components of a bivariate field are given by   $Z_{1}(\bm{x},t)  = a_{11} Y(\bm{x},t,\xi_{1})$ and $Z_{2}(\bm{x},t)  = a_{21} Y(\bm{x},t,\xi_{2}) + a_{22} W(\bm{x},t)$, where $Y(\bm{x},t,\xi)$ is a RF with covariance 
\begin{equation}
\label{KK}
\displaystyle
 \text{K}(\theta,u,\zeta_{ij})  =   \frac{1}{ \left( \frac{400 \theta}{c_{S,1}}  +1\right)^{1/2}   \left(  \frac{\zeta_{ij}}{ \left( \frac{400 \theta}{c_{S,1}} +1\right)}+1\right)^{1/2}   }  
  \displaystyle \exp \left\{  -  \frac{3 |u|/c_{T,1}}{   \left(   \frac{\zeta_{ij}}{ \left( \frac{400 \theta}{c_{S,1}} +1\right)}+1 \right)^{1/2}}      
     \right\},
     \end{equation}
     generated from Equation (\ref{eq_latent2}) and the first entries in Tables \ref{cm} and \ref{bernstein}.
The field $W(\bm{x},t)$ is independent of $Y(\bm{x},t,\xi)$, with covariance function of  Gneiting type as in Equation (\ref{gneiting_porcu}), 
\begin{equation}
\label{RR}
  \text{R}(\theta,u) =  \frac{1}{ \left( \frac{400 \theta}{c_{S,2}} +1\right)^{1/2}     }  
  \exp \left\{  -  \frac{ 3 |u|/ c_{T,2}}{ \left( \frac{400 \theta}{c_{S,2}}  +1 \right)^{1/2}}      
     \right\}. 
     \end{equation}
  Thus,  the covariance function of  $\bm{Z}(\bm{x},t)$, with ${\rm K}$ and ${\rm R}$ defined by (\ref{KK}) and (\ref{RR}), respectively, is given by
  \begin{eqnarray}
\text{C}_{11}(\theta,u)   &  =  &  a_{11}^2\text{K}(\theta,u,0) \nonumber \\
\text{C}_{22}(\theta,u)   & =  &  a_{21}^2 \text{K}(\theta,u,0) + a_{22}^2 \text{R}(\theta,u)  \nonumber\\
\text{C}_{12}(\theta,u)   & = &  a_{11}a_{21}\text{K}(\theta,u,\zeta_{12}).   \label{model_latent1}
\end{eqnarray} 
Note that $\text{R}(\theta,0)< 0.05$ and $\text{R}(0,|u|) < 0.05$,  for $\theta>c_{S,2}$ and $|u|>c_{T,2}$, respectively.   Moreover, since $\text{C}_{12}(0,0) =  a_{11}a_{21}/(\zeta_{12}+1)^{1/2}$, the parameter  $\zeta_{12}$ is related to the cross-scale as well as the collocated correlation coefficient between the fields.

\section{Data Example:  Temperatures and  Precipitations}

\label{sec_data}

We illustrate the use of the proposed models on a climatological data set.  We consider a bivariate space-time data set of surface air  temperatures (Variable 1) and precipitations (Variable 2)  obtained from a climate model provided by NCAR (National Center for Atmospheric Research),  Boulder, CO, USA.    Specifically, the data come from the Community Climate System Model (CCSM4.0) (see \citealp{doi:10.1175/2011JCLI4083.1}).   Temperatures and precipitations    are physically related variables    and  their joint modeling  is crucial in   climatological disciplines.

\begin{figure}
\centering{

\includegraphics[scale=0.22]{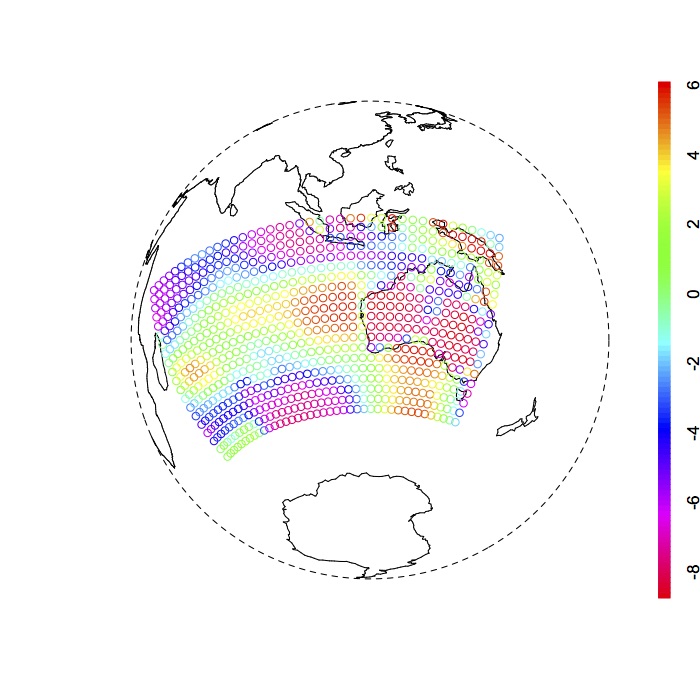} \includegraphics[scale=0.22]{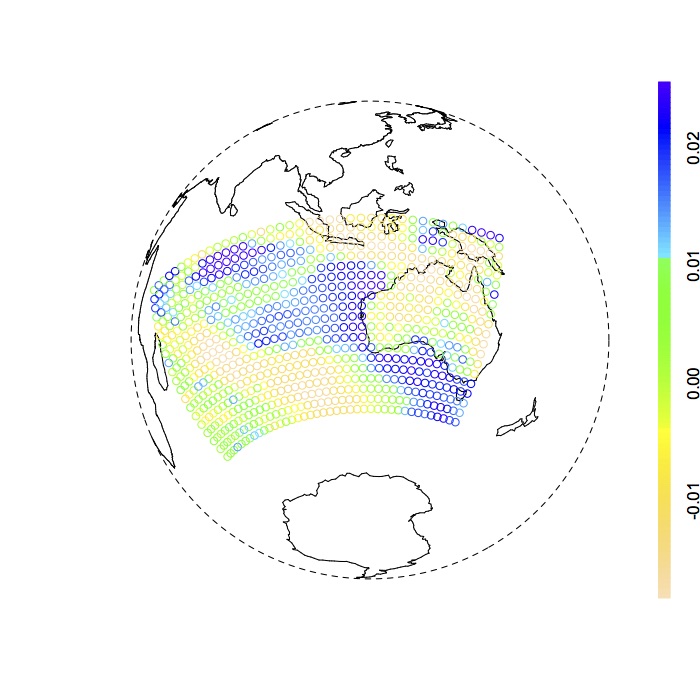}
    }
    \caption{Residuals of the surface air temperatures (left) and precipitations (right) for December 2006, in the area with longitudes between  50 and 150 degrees, and latitudes between $-50$ and 0 degrees.}
\label{mapa}
\end{figure}

In order to attenuate the skewness of precipitations, we work with its cubic root. The units are degrees Kelvin for temperatures and centimeters for the transformed precipitations. The spatial grid   is formed by longitudes and latitudes with $2.5\times 2.5$ degrees of  spatial  resolution (10368 grid points).  We model  planet Earth as a sphere of radius 6378 kilometers.  We focus  on  December, between the years  2006 and 2015 (10 time instants),   and the region with longitudes between  50 and 150 degrees, and latitudes between $-50$ and 0 degrees.  The maximum distance between the points inside this region is approximately half of the maximum distance between any pair of points on planet Earth. Indeed, we are considering a large portion of the globe.   For each variable and  time instant, we use splines to remove the cyclic patterns along longitudes and latitudes. The residuals are approximately Gaussian distributed with zero mean (see Figure \ref{mapa}).    We attenuate the computational burden by taking a sub grid with   10 degrees of resolution in the longitud direction. We thus work with only 200 spatial sites and the resulting data set consists of  4000 observations.    The empirical variances are  $2.033$ and  $4.809\times 10^{-5}$, for Variables 1 and 2, respectively. In addition, these variables are moderately positively correlated. The empirical correlation between the components is $0.279$.

We are interested in showing that in real applications a non-separable model can produce significant improvements with respect to a separable one. For our purposes, we consider the following models:

 \begin{itemize}

\item[(A)] An $m$-separable   covariance  based on the univariate Gneiting class (\ref{gneiting_porcu}):   $$\text{C}_{ij}(\theta,u) = \frac{ \sigma_{ii}\sigma_{jj} \rho_{ij}}{\left( \frac{400\theta}{c_S} +1 \right)^{1/2}} \exp\left\{ \frac{-3|u|/c_T}{\left(\frac{400\theta}{c_S} +1 \right)^{1/2}} \right\},$$
where the vector of parameters is given by $(\sigma_{11}^2,\sigma_{22}^2,\rho_{12}, c_S,c_T)^\top$.
\item[(B)] An $m$-separable version of the modified Gneiting class (\ref{modified1}) with the parameterization $c_S:= c_{11}=c_{22}=c_{12}$. The vector of parameters is the same as in Model A.
\item[(C)]  A non-separable  version of the modified Gneiting class (\ref{modified1}), with the parsimonious parameterization $c_{12}=\max\{c_{11},c_{22}\}$. The vector of parameters is given by $(\sigma_{11}^2,\sigma_{22}^2,\rho_{12}, c_{11},c_{22},c_T)^\top$.
\item[(D)] The non-separable LMC   (\ref{model_latent1}) based on the latent dimension approach. Here, the vector of parameters is given by $(a_{11},a_{21},a_{22}, c_{S,1},c_{S,2},c_{T,1},c_{T,2},\zeta_{12})^\top$.

\end{itemize}

We estimate the models parameters  using the pairwise composite likelihood (CL) method developed for multivariate RFs  by \cite{Bevilacqua2016}. This method is a likelihood approximation  and offers a trade-off between statistical and computational efficiency.    We  only consider  pairs of observations with  spatial distances  less than $1275.6$ kilometers (equivalent to 0.2  radians on the unit sphere) and temporal distance less than $4$ years.

\begin{figure}
\centering{  

\includegraphics[scale=0.25]{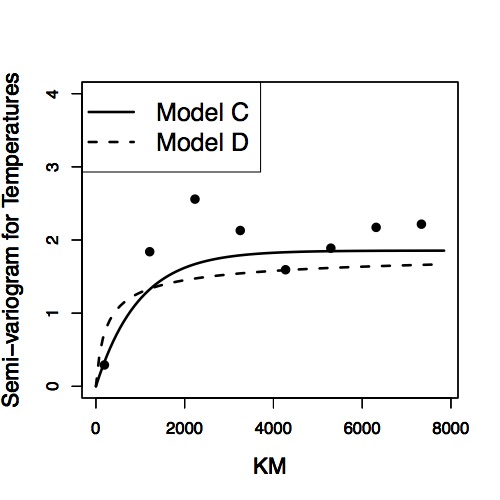}     \includegraphics[scale=0.25]{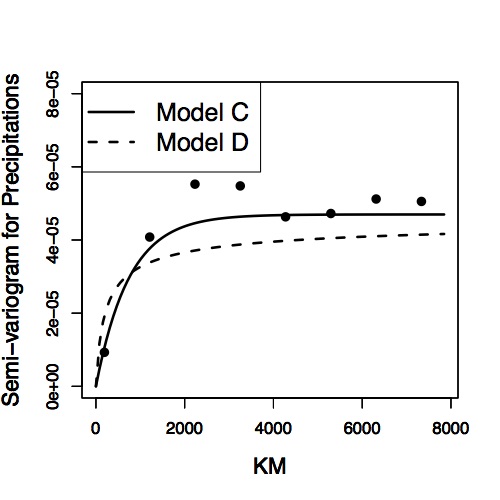}   \includegraphics[scale=0.25]{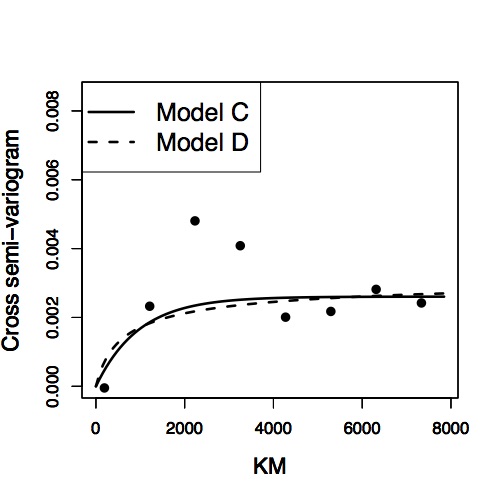}  }

  \includegraphics[scale=0.25]{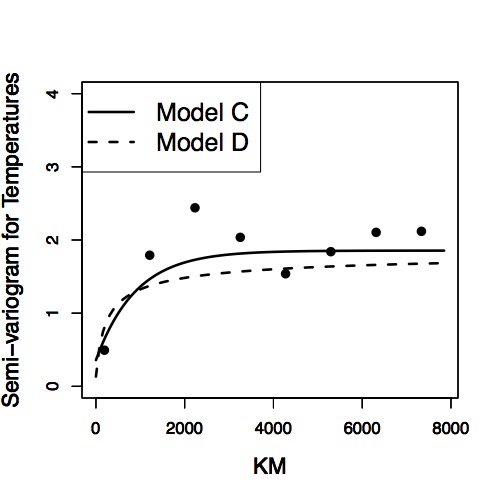}     \includegraphics[scale=0.25]{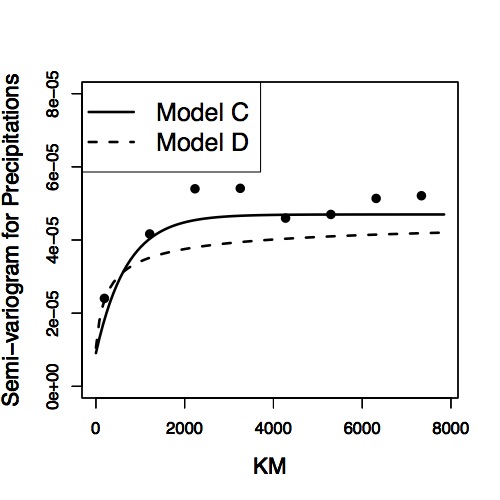}  \includegraphics[scale=0.25]{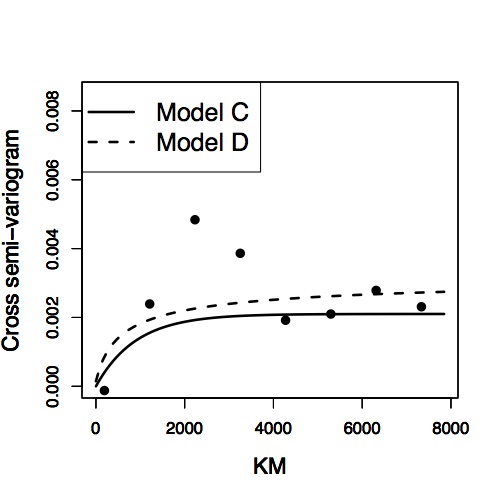}
  \caption{Empirical spatial semi-variograms versus theoretical Models C and D. At the top, we consider the temporal lag as $|u|=0$, whereas at the bottom $|u|=1$.} 
\label{variogram}
\end{figure}

 The CL estimates for   Models A-C   are reported in Table \ref{estimacion_datos2}, whereas  Table \ref{estimacion_datosLMC} reports the results for Model D.    In addition, Table \ref{mspe} contains  the  Log-CL values attained at the optimum  for each model.  As expected, the non-separable covariances C and D exhibit the  highest values of the Log-CL.     In Figure \ref{variogram}, we compare the empirical spatial semi-variograms, for the temporal lags $|u|=0$ and $|u|=1$,  versus the theoretical   models C  and D. Both models seem to capture the behavior of the empirical variograms at the limit values (sills). However, as noted by \cite{Gneiting:Kleibler:Schlather:2010}, disagreements between empirical and theoretical fits are typically   observed in practice, and it can be associated to  biases in the empirical estimators.

 \begin{table}
\centering
\caption{CL estimates  for Models A-C.  For the separable Models A and B, there is a single spatial scale parameter: $c_S:= c_{11}=c_{22}=c_{12}$.  For Model C, the cross scale parameter $c_{12}$ is the maximum between the marginal spatial scales $c_{11}$ and $c_{22}$.}
\label{estimacion_datos2}
\begin{tabular}{ccccccccc} \hline \hline
                        &&  $\widehat{\sigma}_{11}^2$  &  $\widehat{\sigma}_{22}^2$   &  $\widehat{\rho}_{12}$   &  $\widehat{c}_{11}$   &   $\widehat{c}_{22}$    &  $ \widehat{c}_{12}$ & $\widehat{c}_T$     \\          \hline
  Model A       &&    1.85       &   $4.71\times 10^{-5}$  &   0.28       &  102   &   102    &   102   & 11.88         \\ 
  Model B       &&    1.84       &   $4.72\times 10^{-5}$     &     0.28    &       2602   &    2602       &   2602    &   22.58        \\ 
  Model C       &&    1.85       &   $4.69 \times 10^{-5}$    &   0.28    &  2901   &    2245   &   2901  &   22.92     \\  \hline  

\end{tabular}
\end{table}

 \begin{table}
\centering
\caption{CL estimates  for Model D.}
\label{estimacion_datosLMC}
\begin{tabular}{cccccccccc} \hline \hline
                        &&  $\widehat{a}_{11}$     &  $\widehat{a}_{21}$   &  $\widehat{a}_{22}$   &   $\widehat{c}_{S,1}$    &  $\widehat{c}_{S,2}$ & $\widehat{c}_{T,1}$  &  $\widehat{c}_{T,2}$  &  $\widehat{\zeta}_{12}$  \\          \hline
   Model D       &&    1.38       &   $3.98\times 10^{-3}$       &  $5.58\times 10^{-3}$   &    47133   &  30805  &  41.36      &  4.03 & 1.66   \\     
   \hline   
\end{tabular}
\end{table}

Finally, we compare the predictive performance  of the covariances through a cross-validation study based on the cokriging predictor. We use a drop-one prediction strategy and quantify the discrepancy between the real and predicted values for each variable through the   Mean Squared  Error (MSE) and the Continuous Ranked Probability Score (CRPS)  \citep{zhang2010kriging}.  Smaller values of these indicators  imply  better  predictions.     Table \ref{mspe} displays the results, where $\text{MSE}_i$ denotes the $\text{MSE}$ associated to the variable $i$, for $i=1,2$. The interpretation of   $\text{CRPS}_i$ is similar.    Apparently,  the performance of Model A is quite poor in comparison to the other  models.     Moreover, the non-separable  models, C and D,  have better results than the $m$-separable model B.    Note that although Model B has a smaller $\text{MSE}_{2}$ value than Model D,  D is globally superior.      If we choose C instead of B, we have a significant   improvement   in the  prediction of Variable 1, since the corresponding ${\rm MSE}_{1}$ ratio is approximately 0.84. Similarly,   the ${\rm CRPS}_1$ ratio  is 0.93.   Also, direct inspection of  ${\rm MSE}_{2}$ and ${\rm CRPS}_2$ show that  Model C provides improvements in the  prediction of Variable 2 in comparison to Model B.   Indeed, for this specific data set, Model C shows the best results in terms of MSE,  outperforming even Model D.  A nice property of the multivariate modified Gneiting class C  is that it has physically interpretable parameters. Naturally, the proposed models can be combined to  provide more flexible models.

\begin{table}
\centering
\caption{Comparison between Models A-D in terms of   Log-CL and cross-validation scores. Here, ${\rm MSE}_{i}$ and ${\rm CRPS}_i$ denote respectively the MSE and CRPS associated to Variable $i$, for $i=1,2$.}
\label{mspe}
\begin{tabular}{cccccccccc} \hline \hline
                         &    Log-CL $\times 10^{-5}$        &     ${\rm MSE_1}$    &     ${\rm MSE_2}  \times 10^{5}$   &    ${\rm CRPS_1}$   & ${\rm CRPS_2}$       \\ \hline
  Model A  &    5.79   &  0.27      &      2.07     &  1.21   &   1.52        \\
  Model B    &   5.80   &  0.13      &     1.62   &    0.67   &  1.02                   \\
  Model C  &    5.81    &  0.11      &     1.61   &   0.62  &    0.93        \\  
  Model D  &   5.81   &  0.11      &      1.84      &   0.63   &  0.62     \\
\hline 
\end{tabular}
\end{table}

\section{Discussion}
\label{discussion}

In the paper, we  have discussed several  construction principles for multivariate space-time covariances, which  allow to escape from separability. In particular, we have proposed Gneiting type families of cross-covariances and  their properties  have been illustrated through a real data example of surface air temperatures and precipitations.  The proposed models have shown a good performance using separable models as a benchmark. We believe that the methodology used to prove our theoretical results can be adapted to find additional flexible classes of matrix-valued covariances and to develop new applications on  univariate or multivariate global phenomena evolving temporally.  Specifically, our proposals can be used as a building  block for the construction of more sophisticated models, such as non-stationary RFs.

\section*{Acknowledgments}

Alfredo Alegr\'ia's work was initiated when he was a PhD student at   Universidad T\'ecnica Federico Santa Mar\'ia.  Alfredo Alegr\'ia was supported by Beca CONICYT-PCHA/Doctorado Nacional/2016-21160371.  Emilio Porcu is supported by Proyecto Fondecyt Regular number 1170290.   Reinhard Furrer acknowledges support of the Swiss National Science Foundation
 SNSF-144973.  Jorge Mateu is partially supported by grant MTM2016-78917-R. We additionally acknowledge the World Climate Research Programme's Working Group on Coupled modeling, which is responsible for Coupled Model Intercomparison Project (CMIP).

%%%%%%%%%%%%%%%%%%%%%%%%%%%%%%%%%%%%%%%%%%%%%%%%%%%%%%%%%%%%%
%%%%%%%%%%%%%%%%%%%%%%%%%%%%%%%%%%%%%%%%%%%%%%%%%%%%%%%%%%%%%
%%%%%%%%%%%%%%%%%%%%%%%%%%%%%%%%%%%%%%%%%%%%%%%%%%%%%%%%%%%%%
\appendix
\section*{Appendices}
\addcontentsline{toc}{section}{Appendices}
\renewcommand{\thesubsection}{\Alph{subsection}}
\renewcommand{\theequation}{\Alph{subsection}.\arabic{equation}}
\renewcommand\thetable{\thesubsection.\arabic{table}}

\counterwithin{thm}{subsection}
\counterwithin{mydef}{subsection}

\subsection{Background on Positive Definite Functions}

\label{background}

 We start by defining  positive definite  functions, that  arise in statistics as  the  covariances  of Gaussian RFs as well as the characteristic functions of probability distributions.

\begin{mydef}	
Let $\mathcal{E}$ be a non-empty set and $m\in\mathbb{N}$. We say that   the matrix-valued function $\textbf{\text{F}}:\mathcal{E}\times \mathcal{E} \rightarrow \mathbb{R}^{m\times m}$ is positive definite if  for all integer $n\geq 1$, $\{e_1,\hdots,e_n\} \subset \mathcal{E}$ and $\{\bm{a}_1,\hdots,\bm{a}_n\}\subset \mathbb{R}^m$,  the following inequality holds:
\begin{equation} \label{posdef}
   \sum_{\ell=1}^{n} \sum_{r=1}^n   \bm{a}_{\ell}^\top \textbf{\text{F}}(e_\ell,e_r) \bm{a}_{r} \geq 0.
 \end{equation}
We denote as $\mathcal{P}^m(\mathcal{E})$ the class of such mappings $\textbf{\text{F}}$ satisfying Equation (\ref{posdef}).  
\end{mydef}
Next, we focus on the cases where $\mathcal{E}$ is either $\mathbb{R}^d$, $\mathbb{S}^d$ or $\mathbb{S}^d \times \mathbb{R}^k$, for $d,k\in\mathbb{N}$. 
 For a clear presentation of the results, Table \ref{notacion} summarizes the notation introduced along this Appendix.

\subsubsection{Matrix-valued Positive Definite Functions on Euclidean Spaces: The Classes $\Phi^{m}_{d,\mathcal{S}}$ and $\Phi^m_{d,\mathcal{I}}$}

This section provides a brief  review about matrix-valued positive definite  functions  on  the Euclidean space $\mathcal{E} = \mathbb{R}^d$.  Specifically, we expose some characterizations for the stationary and Euclidean isotropic members of the class $\mathcal{P}^m(\mathbb{R}^d)$.

\begin{table}
\begin{center}
\caption{Summary of the notation used along the Appendices. Throughout, in the univariate case ($m=1$) we omit the super index:  $\mathcal{P}(\mathcal{E})$, $\Phi_{d,\mathcal{S}}$, $\Phi_{d,\mathcal{I}}$, $\Psi_{d,\mathcal{I}}$ and $\Upsilon_{d,k}$ are used instead of $\mathcal{P}^1(\mathcal{E})$, $\Phi_{d,\mathcal{S}}^1$, $\Phi_{d,\mathcal{I}}^1$, $\Psi_{d,\mathcal{I}}^1$ and $\Upsilon_{d,k}^1$, respectively. }
\label{notacion}
\begin{tabular}{c|c}\hline \hline
Notation &  Description\\ \hline
$\mathcal{P}^m(\mathcal{E})$ & Positive definite matrix-valued ($m\times m$) functions on $\mathcal{E} \times \mathcal{E}$.  \\  \hline
$\Phi^m_{d,\mathcal{S}}$ &  Continuous, bounded and stationary  elements of $\mathcal{P}^m(\mathbb{R}^d)$.\\  \hline
$\Phi_{d,\mathcal{I}}^m$ &  Continuous, bounded   and Euclidean isotropic  elements of $\mathcal{P}^m(\mathbb{R}^d)$.\\ \hline
$\Psi_{d,\mathcal{I}}^m$ &  Continuous, bounded and geodesically isotropic elements of $\mathcal{P}^m(\mathbb{S}^d)$.\\   \hline
$\Upsilon_{d,k}^m$         &  Elements in $\mathcal{P}^m(\mathbb{S}^d \times \mathbb{R}^k)$ being, continuous, bounded, geodesically \\  
    & isotropic in the     spherical variable    and stationary in the Euclidean variable.\\
  \hline 
\end{tabular}
\end{center}
\end{table}

 We say that $\textbf{\text{F}} \in \mathcal{P}^m(\mathbb{R}^d)$ is stationary if there exists a mapping $\tilde{\bm{\varphi}}:\mathbb{R}^d \rightarrow \mathbb{R}^{m\times m}$ such that 
\begin{equation} \label{eucl_stationary}
\textbf{\text{F}}(\bm{x},\bm{y}) = \tilde{\bm{\varphi}}(\bm{x}-\bm{y}) = [\tilde{\varphi}_{ij}(\bm{x}-\bm{y})]_{i,j=1}^m, \qquad \bm{x},\bm{y}\in\mathbb{R}^d.
\end{equation}
   We call $\Phi_{d,\mathcal{S}}^m$ the class of continuous mappings  $\tilde{\bm{\varphi}}$ such that $\textbf{\text{F}}$ in (\ref{eucl_stationary}) is positive  definite.    Cram\'er's Theorem  \citep{cramer1940theory}  establishes that  $\tilde{\bm{\varphi}}\in \Phi_{d,\mathcal{S}}^m$   if and only if it can be represented through
\begin{equation}\label{cramer}
  \tilde{\bm{\varphi}}(\bm{h}) = \int_{\mathbb{R}^d}  \exp\{ \imath   \bm{h}^\top \bm{\omega}  \} \text{d}\tilde{\bm{\Lambda}}_d(\bm{\omega}),   \qquad \bm{h}\in\mathbb{R}^d,
  \end{equation}
where $ \imath   = \sqrt{-1} \in\mathbb{C}$ and $\tilde{\bm{\Lambda}}_d: \mathbb{R}^d\rightarrow \mathbb{C}^{m\times m}$ is a matrix-valued mapping, with increments being  Hermitian and  positive definite matrices,   and whose elements, $\tilde{\Lambda}_{d,ij}(\cdot)$, for $i,j=1,\hdots,m$, are functions of bounded variation (see \citealp{Wackernagel:2003}).   In particular, the diagonal terms, $\tilde{\Lambda}_{d,ii}(\bm{\omega})$, are real, non-decreasing and bounded, whereas  the off-diagonal elements are generally complex-valued.   Cramer's Theorem  is the multivariate version  of the celebrated Bochner's Theorem \citep{bochner1955harmonic}. If the elements of $\tilde{\bm{\Lambda}}_d(\cdot)$ are absolutely continuous, then Equation (\ref{cramer}) simplifies to
\begin{equation*}\label{cramer2}
  \tilde{\bm{\varphi}}(\bm{h}) = \int_{\mathbb{R}^d}  \exp\{ \imath   \bm{h}^\top \bm{\omega}  \}  \tilde{\bm{\lambda}}_d(\bm{\omega})\text{d}\bm{\omega},   \qquad \bm{h}\in\mathbb{R}^d,
  \end{equation*}
with  $\tilde{\bm{\lambda}}_d(\bm{\omega})=[\tilde{\lambda}_{d,ij}(\bm{\omega})]_{i,j=1}^m$ being  Hermitian and positive definite, for any $\bm{\omega}\in\mathbb{R}^d$. The mapping $\tilde{\bm{\lambda}}_d(\bm{\omega})$ is known as the matrix-valued  \textit{spectral density}  and classical Fourier inversion yields 
$$   \tilde{\bm{\lambda}}_d(\bm{\omega}) = \frac{1}{(2\pi)^d} \int_{\mathbb{R}^d}  \exp\{ - \imath   \bm{h}^\top \bm{\omega}  \}  \tilde{\bm{\varphi}}(\bm{h})  \text{d}\bm{h}, \qquad \bm{\omega}\in\mathbb{R}^d.  $$

Finally, the following inequality between the elements of $\tilde{\bm{\varphi}}$ is true
$$   |\tilde{\varphi}_{ij}(\bm{h})|^2 \leq \tilde{\varphi}_{ii}(\bm{0}) \tilde{\varphi}_{jj}(\bm{0}),\qquad  \text{ for all } \bm{h}\in\mathbb{R}^d. $$
However, the maximum value of the mapping $\tilde{\varphi}_{ij}(\bm{h})$, with $i\neq j$, is not necessarily reached at  $\bm{h}=\bm{0}$.  In general, $\tilde{\varphi}_{ij}$ is not itself a scalar-valued positive definite function when $i\neq j$.

 Consider an element $\textbf{\text{F}}$ in $\mathcal{P}^m(\mathbb{R}^d)$ and suppose that there exists a continuous and bounded mapping $\bm{\varphi}:\mathbb{R}_+ \rightarrow \mathbb{R}^{m\times m}$ such that  
\begin{equation*}\label{eucl_isotropic}
\textbf{\text{F}}(\bm{x},\bm{y}) = \bm{\varphi}(\|\bm{x}-\bm{y}\|), \qquad \bm{x},\bm{y}\in\mathbb{R}^d.
\end{equation*}
 Then, $\textbf{\text{F}}$ is called stationary and \textit{Euclidean isotropic} (or radial). We denote as $\Phi^m_{d,\mathcal{I}}$ the class of  bounded, continuous, stationary and Euclidean isotropic mappings  $\bm{\varphi}(\cdot)=[\varphi_{ij}(\cdot)]_{i,j=1}^m$.    
 
 When $m=1$, characterization of the class $\Phi_{d,\mathcal{I}}$ was provided through the  celebrated paper by \cite{10.2307/1968466}. 
\cite{AlonsoMalaver2015251} characterize the class $\Phi^m_{d,\mathcal{I}}$ through the continuous members  $\bm{\varphi}$ having representation
\begin{equation*}\label{schoenberg_euclid}
  \bm{\varphi}(r) = \int_{[0,\infty)}  \Omega_d(r\omega) \text{d}\bm{\Lambda}_d(\omega),  \qquad  r\geq 0,
  \end{equation*}
where $\bm{\Lambda}_d: [0,\infty)\rightarrow \mathbb{R}^{m\times m}$ is a matrix-valued mapping, with  increments being positive definite matrices, and  elements $\Lambda_{d,ij}(\cdot)$  of bounded variation, for each  $i,j=1,\hdots,m$.  Here, the function $\Omega_d(\cdot)$ is defined as
 \begin{equation*}
  \label{hankel} 
  \Omega_d(z) = \Gamma(d/2)(z/2)^{-(d-2)/2} J_{(d-2)/2}(z), \qquad z \ge 0,
   \end{equation*}
with  $\Gamma$ being the Gamma function and $J_\nu$  the Bessel function of the first kind of degree $\nu$ (see \citealp{Abramowitz-Stegun:1965}).   If the elements of $\bm{\Lambda}_d(\cdot)$ are absolutely continuous, then we have an associated spectral density $\bm{\lambda}_d:[0,\infty)\rightarrow \mathbb{R}^{m\times m}$ as in the stationary case, which is called, following \cite{daley2014dimension}, a $d$-Schoenberg matrix.

The classes $\Phi^m_{d,\mathcal{I}}$ are non-increasing in $d$, and  the following inclusion relations are strict
\begin{equation*}\label{inclusion_eucl}
\Phi^m_{\infty,\mathcal{I}} := \bigcap_{d=1}^\infty  \Phi^m_{d,\mathcal{I}}  \subset \cdots \subset \Phi^m_{2,\mathcal{I}} \subset \Phi^m_{1,\mathcal{I}}.
\end{equation*}
The elements $\bm{\varphi}$ \ in the class $\Phi^m_{\infty,\mathcal{I}}$ can be represented as 
\begin{equation*}
  \bm{\varphi}(r) = \int_{[0,\infty)}  \exp(-r^2\omega^2) \text{d}\bm{\Lambda}(\omega),  \qquad  r\geq 0,
  \end{equation*}
where $\bm{\Lambda}$ is a matrix-valued mapping with similar properties as $\bm{\Lambda}_d$.

\subsubsection{Matrix-valued Positive Definite Functions on $\mathbb{S}^d$:  The Class $\Psi_{d,\mathcal{I}}^m$}

In this section, we pay attention to matrix-valued  positive definite functions on the unit sphere. Consider $\textbf{\text{F}} = [\text{F}_{ij}]_{i,j=1}^m\in \mathcal{P}^m(\mathbb{S}^d)$.   We say that $\textbf{\text{F}}$  is \textit{geodesically isotropic}   if there exists a bounded and continuous mapping  $\bm{\psi}:[0,\pi]\rightarrow \mathbb{R}^{m\times m}$ such that
\begin{equation*}\label{sph_isotropic}
\textbf{\text{F}}(\bm{x},\bm{y}) = \bm{\psi}(\theta(\bm{x},\bm{y})), \qquad \bm{x},\bm{y}\in\mathbb{S}^d.
\end{equation*}
The continuous mappings $\bm{\psi}$  are the elements of the class $\Psi_{d,\mathcal{I}}^m$ and  the following inclusion relations are true:
 \begin{equation}\label{inclusion_sph}
  \Psi_{\infty,\mathcal{I}}^m  = \bigcap_{d= 1}^\infty \Psi_{d,\mathcal{I}}^m  \subset \cdots \subset \Psi_{2,\mathcal{I}}^m \subset  \Psi_{1,\mathcal{I}}^m,
  \end{equation}
 where $\Psi_{\infty,\mathcal{I}}^m$ is the class of  geodesically isotropic positive definite functions being valid on the Hilbert sphere $\mathbb{S}^\infty = \{ (x_n)_{n\in\mathbb{N}} \in \mathbb{R}^{\mathbb{N}} : \sum_{n\in\mathbb{N}} x_n^2 =1 \}$.

The elements of the class $\Psi_{d,\mathcal{I}}^m$ have an explicit connection with Gegenbauer (or ultraspherical) polynomials \citep{Abramowitz-Stegun:1965}. 
 Here,  $\mathcal{G}_n^{\lambda}$ denotes  the $\lambda$-Gegenbauer  polynomial of degree $n$, which is defined implicitly through the expression
 $$  \frac{1}{(1+r^2-2r\cos\theta)^\lambda} = \sum_{n=0}^\infty  r^n \mathcal{G}_n^{\lambda}(\cos\theta), \qquad \theta\in[0,\pi], \qquad r\in(-1,1).$$
In particular, $\mathcal{T}_n := \mathcal{G}_n^{0}$ and $\mathcal{P}_n := \mathcal{G}_n^{1/2}$ are respectively the Chebyshev and Legendre polynomials of degree $n$.

 The following result (\citealp{hannan2009multiple}; \citealp{Yaglom:1987})  offers a complete characterization of the classes $\Psi_{d,\mathcal{I}}^m$ and $\Psi_{\infty,\mathcal{I}}^m$, and
 corresponds to the multivariate version of  Schoenberg's Theorem  \citep{schoenberg1942}. Equalities and summability conditions for matrices must be understood  in a componentwise sense.
 
 \begin{thm}
\label{hannan}
 Let $d$ and $m$ be positive integers.
 \begin{itemize}
\item[(1)]   The mapping  $\bm{\psi}$ is a member  of the class $\Psi_{d,\mathcal{I}}^m$ if and only if it admits the representation
\begin{equation}\label{sch_representation}
\bm{\psi}(\theta) = \sum_{n=0}^\infty \textbf{\text{B}}_{n,d} \frac{\mathcal{G}_n^{(d-1)/2}(\cos\theta)}{\mathcal{G}_n^{(d-1)/2}(1)}, \qquad \theta \in [0,\pi],
\end{equation}
where $\{\textbf{\text{B}}_{n,d}\}_{n=0}^\infty$ is a sequence of symmetric, positive definite and summable  matrices.

\item[(2)]  The mapping $\bm{\psi}$ is a member  of the class $\Psi_{\infty,\mathcal{I}}^m$ if and only if it can be represented as
\begin{equation*}\label{sch_representation_inf}
\bm{\psi}(\theta) = \sum_{n=0}^\infty \textbf{\text{B}}_{n} (\cos\theta)^n, \qquad \theta \in [0,\pi],
\end{equation*}
where $\{\textbf{\text{B}}_{n}\}_{n=0}^\infty$ is a sequence of symmetric, positive definite and summable  matrices.
\end{itemize} 
\end{thm}

Using orthogonality properties of Gegenbauer polynomials \citep{Abramowitz-Stegun:1965} and through classical Fourier inversion we can prove that
\begin{eqnarray}
\label{sch_matrices}
\textbf{\text{B}}_{0,1} & = & \frac{1}{\pi}\int_0^\pi \bm{\psi}(\theta)\text{d}\theta,  \nonumber \\
\textbf{\text{B}}_{n,1} & = & \frac{2}{\pi}\int_0^\pi \cos(n\theta)\bm{\psi}(\theta)\text{d}\theta,  \qquad  \text{ for } n\geq 1,
\end{eqnarray}
whereas  for $d\geq 2$, we have
\begin{equation}
\label{sch_matrices2}
\textbf{\text{B}}_{n,d} = \frac{(2n+d-1)}{2^{3-d}\pi} \frac{[\Gamma((d-1)/2)]^2}{\Gamma(d-1)} \int_{0}^\pi  \mathcal{G}_{n}^{(d-1)/2}(\cos \theta) (\sin\theta)^{d-1} \bm{\psi}(\theta) \text{d}\theta,  \qquad n\geq 0,
\end{equation}
where integration is taken componentwise. The matrices $\{\textbf{\text{B}}_{n,d}\}_{n=0}^\infty$ are called  \textit{Schoenberg's matrices}. For the case $m=1$, such result is reported by \cite{gneiting2013}.

\subsubsection{Matrix-valued Positive Definite Functions on $\mathbb{S}^d\times \mathbb{R}^k$: The Class $\Upsilon_{d,k}^m$}

Let $d$, $k$ and $m$ be  positive integers. We now focus on the class of matrix-valued positive definite functions on  $\mathbb{S}^d \times \mathbb{R}^k$, being bounded,   continuous,    geodesically isotropic in the spherical component  and stationary in the Euclidean one. The case $k=1$  is particularly important, since $\mathcal{P}^m(\mathbb{S}^d\times \mathbb{R})$   can be  interpreted as the class of admissible  space-time covariances for   multivariate Gaussian RFs, with spatial locations on the unit sphere.

Consider $\textbf{\text{F}} \in \mathcal{P}^m(\mathbb{S}^d\times \mathbb{R}^k)$ and suppose that  there exists a bounded and continuous mapping $\textbf{\text{C}}:[0,\pi]\times \mathbb{R}^k \rightarrow \mathbb{R}^{m\times m}$ such that 
\begin{equation}\label{product_space2}
  \textbf{\text{F}}((\bm{x},\bm{t}),(\bm{y},\bm{s})) = \textbf{\text{C}}(\theta(\bm{x},\bm{y}),\bm{t}-\bm{s}),         \qquad \bm{x},\bm{y}\in\mathbb{S}^d, \bm{t},\bm{s}\in\mathbb{R}^k. 
\end{equation}
Such mappings $\textbf{\text{C}}$  are the elements of the class $\Upsilon_{d,k}^m$.   These classes are non-increasing in $d$ and we have the inclusions
   $$  \Upsilon_{\infty,k}^m := \bigcap_{d=1}^\infty   \Upsilon_{d,k}^m  \subset \cdots \subset \Upsilon_{2,k}^m  \subset  \Upsilon_{1,k}^m.$$

\cite{Ma2016} proposes the generalization of Theorem \ref{hannan}  to the space-time case.  Theorem \ref{characterization} below offers a complete characterization of the class $\Upsilon_{d,k}^m$ and $\Upsilon_{\infty,k}^m$, for any $m\geq 1$.  Again, equalities and summability conditions must be understood in a componentwise sense.

\begin{thm} \label{characterization}
Let $d$, $k$ and $m$ be  positive integers and $\textbf{\text{C}}: [0,\pi] \times \R^k \to \mathbb{R}^{m\times m}$  a continuous matrix-valued  mapping, with $\text{C}_{ii}(0,\bm{0})<\infty$, for all $i=1,\hdots,m$.
\begin{itemize}
\item[(1)] The mapping $\textbf{\text{C}}$ belongs to the class $\Upsilon_{d,k}^m$ if and only if 
\begin{equation} \label{2.1}
\textbf{\text{C}}(\theta,\bm{u}) = \sum_{n=0}^{\infty} \tilde{\bm{\varphi}}_{n,d}(\bm{u}) \frac{\mathcal{G}_n^{(d-1)/2}(\cos\theta)}{\mathcal{G}_n^{(d-1)/2}(1)}, \qquad (\theta,\bm{u}) \in [0,\pi] \times \R^k,
\end{equation} 
with $\left \{ \tilde{\bm{\varphi}}_{n,d}(\cdot) \right\}_{n=0}^{\infty}$ being a sequence of members of the class $\Phi_{k,\mathcal{S}}^m$, with the additional requirement that the sequence of matrices $\{ \tilde{\bm{\varphi}}_{n,d}(\bm{0})\}_{n=0}^\infty$ is summable. 
 
 \item[(2)] The mapping $\textbf{\text{C}}$ belongs to the class $\Upsilon_{\infty,k}^m$ if and only if 
\begin{equation} \label{2.2}
\textbf{\text{C}}(\theta,\bm{u}) = \sum_{n=0}^{\infty} \tilde{\bm{\varphi}}_n(\bm{u}) (\cos\theta)^n, \qquad (\theta,\bm{u}) \in [0,\pi] \times \R^k,
\end{equation} 
with $\left \{ \tilde{\bm{\varphi}}_n(\cdot) \right\}_{n=0}^{\infty}$ being a sequence of members of the class $\Phi_{k,\mathcal{S}}^m$, with the additional requirement that the sequence of matrices $\{\tilde{\bm{\varphi}}_n(\bm{0})\}_{n=0}^\infty$ is summable.
 
 \end{itemize}
\end{thm}

Again, using    orthogonality arguments,  we have 
\begin{eqnarray}
\label{sch_functions}
\tilde{\bm{\varphi}}_{0,1}(\bm{u}) & = & \frac{1}{\pi}\int_0^\pi \textbf{\text{C}}(\theta,\bm{u}) \text{d}\theta,  \nonumber \\
\tilde{\bm{\varphi}}_{n,1}(\bm{u}) & = & \frac{2}{\pi}\int_0^\pi \cos(n\theta)\textbf{\text{C}}(\theta,\bm{u}) \text{d}\theta,  \qquad  \text{ for } n\geq 1, 
\end{eqnarray}
whereas  for $d\geq 2$,
\begin{equation} 
\label{sch_functions2}
\tilde{\bm{\varphi}}_{n,d}(\bm{u}) = \frac{(2n+d-1)}{2^{3-d}\pi} \frac{[\Gamma((d-1)/2)]^2}{\Gamma(d-1)}  \int_{0}^{\pi} \mathcal{G}_{n}^{(d-1)/2}(\cos \theta) (\sin \theta)^{d-1} \textbf{\text{C}}(\theta,\bm{u})  {\rm d} \theta, \qquad   n\geq 0.
\end{equation}

\subsection{Proof of Theorem \ref{los_tanVores}}
\label{modified_proof}

In order to illustrate the results following subsequently,  a technical Lemma will be useful. We do not provide a proof because it is obtained following the same arguments as in \cite{Porcu20111293}.
\begin{lem} \label{amorio}
Let $m$, $d$ and $k$ be strictly positive integers. Let $(X,{\cal B}, \mu)$ be a measure space, for $X \subset \R$ and  $\cal B$ being the Borel sigma algebra. Let $\bm{\psi}: [0,\pi] \times X \to \R$ and $\bm{\varphi}: [0,\infty) \times X \to \R$ be continuous mappings satisfying 
\begin{enumerate}
\item $\bm{\psi}(\cdot,\xi) \in \Psi_{d,{\cal I}}^m$ a.e. $\xi \in X$; 
\item $\bm{\psi}(\theta,\cdot) \in L_1(X,{\cal B},\mu)$ for any $\theta \in [0,\pi]$;
\item $\bm{\varphi}(\cdot,\xi) \in \Phi_{k,{\cal I}}^m$ a.e. $\xi \in X$; 
\item $\bm{\varphi}(u, \cdot) \in  L_1(X,{\cal B},\mu)$ for any $u \in [0,\infty)$.
\end{enumerate}
Let $\textbf{\text{C}} :[0,\pi] \times [0,\infty)\rightarrow \mathbb{R}^{m\times m}$ be the mapping defined through
\begin{equation} \label{la_traiccion}
\textbf{\text{C}} (\theta,u) = \int_X \bm{\psi}(\theta,\xi) \bm{\varphi}(u, \xi)  \mu ({\rm d} \xi), \qquad (\theta,u) \in  [0,\pi] \times [0,\infty).
\end{equation}
Then, $\textbf{\text{C}} $ is continuous and bounded. Further, $\textbf{\text{C}} $ belongs to the class $\Upsilon_{d,k}^m$.
\end{lem}
Of course, Lemma \ref{amorio} is a particular case of the scale mixtures introduced in Section \ref{preliminares}.

{\bf Proof of Theorem \ref{los_tanVores}} Let $(X, {\cal B}, \mu)$ as in Lemma \ref{amorio} and consider $X = \R_+$ with $\mu$ the Lebesgue measure.  We offer a proof of the constructive type. Let us define the function
$\bm{\psi}(\theta,\xi)$ with members $\psi_{ij}(\cdot,\cdot)$ defined through 
$$ \psi_{ij}(\theta,\xi) = { \sigma_{ii}\sigma_{jj} \rho_{ij}}
 \left ( 1 - \frac{\theta}{\xi c_{ij} } \right )_{+}^{n+1}, \qquad (\theta, \xi) \in [0,\pi] \times X, \quad i,j=1,\ldots,m, $$
where, as asserted, the constants $\sigma_{ii}$, $\rho_{ij}$ and $c_{ij}$ are determined according to condition (\ref{askey3}). 
Let us now define the mapping $(u,\xi) \mapsto \varphi(u,\xi):= \xi^{n+1} (1 - \xi  f(u))_{+}^{\ell}$, with $(u,\xi) \in [0,\infty) \times X$. It can be verified that both $\bm{\psi}$ and $\varphi$ satisfy requirements 1--4 in Lemma \ref{amorio}. In particular, Condition 1 yields thanks to Lemmas 3 and 4 in \cite{gneiting2013}, as well as Theorem 1 in \cite{daley2015}. Also, arguments in \cite{PBG16}  show that Condition 3 holds for any $\ell \ge 1$. We can now apply Lemma \ref{amorio}, so that we have that
$$ \text{C}_{i,j, n,\ell}(\theta,u):= \int_X \bm{\psi}(\theta,\xi) \varphi(u,\xi) {\rm d} \xi, \qquad [0, \pi] \times [0,\infty) $$
is a member of the class $\Upsilon_{2n+1,1}^m$ for any $\ell \ge 1$. Pointwise application of an elegant scale mixture argument as in Proposition 1 of \cite{PBG16}  shows that 
\begin{equation}
 \text{C}_{i,j,n,\ell} (\theta,u) = \mathcal{B}(n+2,\ell+1) \frac{ \rho_{ij} \sigma_{ii} \sigma_{jj}}{f(u)^{n+2}}  \left( 1- \frac{\theta f(u)}{c_{ij}} \right)_+^{n+\ell+1}, \qquad (\theta,u) \in  [0,\pi] \times [0,\infty),
\end{equation}
where $\mathcal{B}$ denotes the Beta function \citep{Abramowitz-Stegun:1965}.  We now omit the factor $\mathcal{B}(n+2,\ell+1)$ since it does not affect positive definiteness. Now, standard convergence arguments show that $$\lim_{\ell \to \infty} \text{C}_{i,j,n,\ell}(\theta / \ell,u) =  \frac{ \rho_{ij} \sigma_{ii} \sigma_{jj}}{f(u)^{n+2}}  \exp \left ( -  \frac{\theta f(u)}{c_{ij}} \right ), (\theta,u) \in [0,\pi] \times [0,\infty), $$ with the convergence being uniform in any compact set. The proof is then completed in view of Bernstein's theorem  \citep{feller1966introduction}.

%\subsection{Proof of Theorem \ref{prop_latent1}}
%\label{latent1}

\subsection{Proof of Theorem \ref{prop_latent2}}
\label{latent2}

Before we state the proof of Theorem \ref{prop_latent2}, we need to introduce three auxiliary lemmas.

\begin{lem}
\label{primer_lemma}
Let $\text{C}: [0,\pi]\times   \mathbb{R}^k \rightarrow \mathbb{R}$ be  a continuous, bounded and integrable function, for some positive integer  $k$. Then $\text{C} \in \Upsilon_{d,k}$, for  $d\geq1$, if and only if the mapping $\psi_{\bm{\omega}}:[0,\pi]  \rightarrow  \mathbb{R}$ defined as 
\begin{equation}
\label{c_w1}
\psi_{\bm{\omega}}(\theta) = \frac{1}{(2\pi)^k}\int_{\mathbb{R}^k} \exp\{-\imath \bm{\omega}^\top \bm{v} \}  \text{C}(\theta,\bm{v}) \text{d}\bm{v}, \qquad \theta\in [0,\pi] ,
\end{equation}
belongs to the class  $\Psi_{d,\mathcal{I}}$, for all $\bm{\omega}\in\mathbb{R}^k$.
\end{lem}
We do not report the proof of Lemma \ref{primer_lemma} since the arguments are the same as in the proof of Lemma \ref{lema_2} below. Note that this lemma is a spherical version of the  result given by \cite{Cressie:Huang:1999}.

\begin{lem}
\label{lema_2}
Let $\text{C}: [0,\pi]\times  \mathbb{R}^l \times \mathbb{R}^k \rightarrow \mathbb{R}$ be  a continuous, bounded and integrable function, for some positive integers  $l$ and $k$. Then, $\text{C} \in \Upsilon_{d,k+l}$, with $d\geq1$, if and only if the mapping $\text{C}_{\bm{\omega}}:[0,\pi] \times \mathbb{R}^l  \rightarrow  \mathbb{R}$ defined as 
\begin{equation}
\label{c_w}
\text{C}_{\bm{\omega}}(\theta,\bm{u}) = \frac{1}{(2\pi)^k}\int_{\mathbb{R}^k} \exp\{-\imath \bm{\omega}^\top \bm{v} \}  \text{C}(\theta,\bm{u},\bm{v}) \text{d}\bm{v}, \qquad (\theta,\bm{u})\in [0,\pi] \times \mathbb{R}^l,
\end{equation}
belongs to the class  $\Upsilon_{d,l}$, for all $\bm{\omega}\in\mathbb{R}^k$.
\end{lem}

\textbf{Proof of Lemma \ref{lema_2}} Suppose that $\text{C} \in \Upsilon_{d,k+l}$, then the characterization of \cite{berg2016schoenberg} implies that 
\begin{equation}
\text{C}(\theta,\bm{u},\bm{v}) =  \sum_{n=0}^\infty \tilde{\varphi}_{n,d}(\bm{u},\bm{v}) \frac{\mathcal{G}_n^{(d-1)/2}(\cos\theta)}{\mathcal{G}_n^{(d-1)/2}(1)},
\end{equation}
where $\{\tilde{\varphi}_{n,d}(\cdot,\cdot)\}_{n=0}^\infty$ is a sequence of functions in $\Phi_{k+l,\mathcal{S}}$, with $ \sum_{n=0}^\infty \tilde{\varphi}_{n,d}(\bm{0},\bm{0}) < \infty$.
Therefore,
\begin{eqnarray*}
\text{C}_{\bm{\omega}}(\theta,\bm{u}) &  = &  \frac{1}{(2\pi)^k} \int_{\mathbb{R}^k} \exp\{-\imath \bm{\omega}^\top \bm{v} \}  \left( \sum_{n=0}^\infty \tilde{\varphi}_{n,d}(\bm{u},\bm{v}) \frac{\mathcal{G}_n^{(d-1)/2}(\cos\theta)}{\mathcal{G}_n^{(d-1)/2}(1)} \right)  \text{d}\bm{v}\\
   &   =   &     \sum_{n=0}^\infty     \left(  \frac{1}{(2\pi)^k}\int_{\mathbb{R}^k} \exp\{-\imath \bm{\omega}^\top \bm{v} \}  \tilde{\varphi}_{n,d}(\bm{u},\bm{v})    \text{d}\bm{v}   \right)      \frac{\mathcal{G}_n^{(d-1)/2}(\cos\theta)}{\mathcal{G}_n^{(d-1)/2}(1)},  
 \end{eqnarray*}
where the last step is justified by dominated convergence.  We need to prove that for each fixed $\bm{\omega}\in\mathbb{R}^k$, the sequence of  functions 
$$  \bm{u} \mapsto \tilde{\lambda}_{n,d}(\bm{u};\bm{\omega}) : = \frac{1}{(2\pi)^k} \int_{\mathbb{R}^k} \exp\{-\imath \bm{\omega}^\top \bm{v} \}  \tilde{\varphi}_{n,d}(\bm{u},\bm{v})    \text{d}\bm{v}, \qquad n\geq 0, $$
belongs to the class $\Phi_{l,\mathcal{S}}$, a.e. $\bm{\omega}\in\mathbb{R}^k$. In fact, we have that
\begin{equation}
\label{bochner}
   \frac{1}{(2\pi)^{l}}  \int_{\mathbb{R}^l}   \exp\{-\imath \bm{\tau}^\top \bm{u} \}  \tilde{\lambda}_{n,d}(\bm{u};\bm{\omega}) \text{d}\bm{u}  =   \frac{1}{(2\pi)^{k+l}}  \int_{\mathbb{R}^l}  \int_{\mathbb{R}^k} \exp\{-\imath \bm{\tau}^\top \bm{u} - \imath  \bm{\omega}^\top \bm{v} \}  \tilde{\varphi}_{n,d}(\bm{u},\bm{v})    \text{d}\bm{v}  \text{d}\bm{u}.
   \end{equation}
Since $\tilde{\varphi}_{n,d}(\cdot,\cdot)$  belongs to $\Phi_{k+l,\mathcal{S}}$,  Bochner's Theorem implies that  the right side in Equation (\ref{bochner}) is non-negative everywhere.  This implies that $\tilde{\lambda}_{n,d}(\cdot;\bm{\omega})$ belongs to  the class  $\Phi_{l,\mathcal{S}}$. Also, direct inspection shows that $\sum_{n=0}^\infty \tilde{\lambda}_{n,d}(\bm{0};\bm{\omega}) < \infty$, for all $\bm{\omega}\in\mathbb{R}^k$. The necessary part is completed.

On the other hand, suppose that for each $\bm{\omega}\in\mathbb{R}^k$ the function $\text{C}_{\bm{\omega}}(\theta,\bm{u})$ belongs to the class $\Upsilon_{d,l}$, then there exists a sequence of mappings $\{ \tilde{\lambda}_{n,d}(\cdot;\bm{\omega}) \}_{n=0}^\infty$ in $\Phi_{l,\mathcal{S}}$  for each $\bm{\omega}\in\mathbb{R}^k$, such that
\begin{equation*}
\text{C}_{\bm{\omega}}(\theta,\bm{u}) =  \sum_{n=0}^\infty \tilde{\lambda}_{n,d}(\bm{u};\bm{\omega})  \frac{\mathcal{G}_n^{(d-1)/2}(\cos\theta)}{\mathcal{G}_n^{(d-1)/2}(1)}.
\end{equation*}
Thus,
\begin{eqnarray*}
\text{C}(\theta,\bm{u},\bm{v})    &  = &  \int_{\mathbb{R}^k} \exp\{ \imath \bm{\omega}^\top \bm{v} \}  \left(  \sum_{n=0}^\infty \tilde{\lambda}_{n,d}(\bm{u};\bm{\omega})   \frac{\mathcal{G}_n^{(d-1)/2}(\cos\theta)}{\mathcal{G}_n^{(d-1)/2}(1)} \right)  \text{d}\bm{\omega}\\
&  =  &    \sum_{n=0}^\infty     \left(    \int_{\mathbb{R}^k} \exp\{ \imath \bm{\omega}^\top \bm{v} \} \tilde{\lambda}_{n,d}(\bm{u};\bm{\omega})  \text{d}\bm{\omega} \right)        \frac{\mathcal{G}_n^{(d-1)/2}(\cos\theta)}{\mathcal{G}_n^{(d-1)/2}(1)}.
 \end{eqnarray*}
 We conclude the proof by invoking again  Bochner's Theorem and the result of \cite{berg2016schoenberg}.

\begin{lem}
\label{prop_latent1}
Let $d$ and $k$ be two positive integers. Consider  $g$ and $f$ be  completely monotone and  Bernstein functions, respectively. Then,
\begin{equation}
\text{K}(\theta,\bm{v}) =  \frac{1}{    \left\{  f(\theta) |_{[0,\pi]}  \right\}^{k/2}      } g\left(   \frac{\|\bm{v}\|^2}{   f(\theta) |_{[0,\pi]}    } \right), \qquad  (\theta,\bm{v}) \in [0,\pi]\times\mathbb{R}^k,
\end{equation}
belongs to the class $\Upsilon_{d,k}$, for any positive integer $d$.
\end{lem}

\textbf{Proof of Lemma \ref{prop_latent1}} By Lemma \ref{primer_lemma}, we must  show that $\psi_{\bm{\omega}}$,  defined through Equation (\ref{c_w1}), belongs to the class  $\Psi_{d,\mathcal{I}}$, for all $\bm{\omega}\in\mathbb{R}^k$.  In fact,  we can assume that $\text{C}$ is integrable, since the general case is obtained with the same arguments given by \cite{doi:10.1198/016214502760047113}. Bernstein's Theorem establishes that $g$ can be represented as the Laplace transform of a  bounded measure $G$, then 
\begin{eqnarray*}
                      \psi_{\bm{\omega}} (\theta)   &  =  &      \int_{\mathbb{R}^k} \exp\{-\imath \bm{\omega}^\top \bm{v} \}      \frac{1}{    \left\{  f(\theta) |_{[0,\pi]}  \right\}^{k/2}      }   \int_{[0,\infty)}  \exp\left\{ -   \frac{r\|\bm{v}\|^2}{   f(\theta) |_{[0,\pi]}    } \right\} \text{d}G(r) \text{d}\bm{v}\\
                       &  =  &  \pi^{k/2}    \int_{[0,\infty)}   \exp\left\{ -  \frac{\|\bm{\omega}\|^2}{4r}  f(\theta)  |_{[0,\pi]}  \right\} \text{d}\tilde{G}(r),
\end{eqnarray*}
where   the last equality follows from Fubini's Theorem and $\text{d}G(r) = r^{k/2} \text{d}\tilde{G}(r)$.  In addition, the composition between a negative exponential and a Bernstein function is completely monotone on the real line \citep{feller1966introduction}. Then, for any $\bm{\omega}$ and $r$, the mapping $\theta \mapsto  \exp\{ -  \|\bm{\omega}\|^2  f(\theta)|_{[0,\pi]} /(4r)\}$ is the restriction of a completely monotone function to the interval $[0,\pi]$. Theorem 7 in \cite{gneiting2013} implies that such mapping,  and thus $ \psi_{\bm{\omega}}$, belongs to the class $\Psi_{d,\mathcal{I}}$, for any $d\in\mathbb{N}$ and $\bm{\omega}\in\mathbb{R}^k$.

\textbf{Proof of Theorem \ref{prop_latent2}} By Lemma \ref{lema_2}, we must  show that (\ref{c_w}) belongs to the class $\Upsilon_{d,l}$, for all $d\in\mathbb{N}$. In fact,   assuming again that $\text{C}$ is integrable and invoking  Bernstein's Theorem we have
\begin{eqnarray*}
                       \text{C}_{\bm{\omega}} (\theta,\bm{u})   &  =  &      \int_{\mathbb{R}^k} \exp\{-\imath \bm{\omega}^\top \bm{v} \}   \frac{1}{  \{f_2( \theta )|_{[0,\pi]}\}^{l/2} \left\{f_1\left[  \frac{\|\bm{u}\|^2}{  f_2( \theta )|_{[0,\pi]} } \right]\right\}^{k/2}}   \int_{[0,\infty)}  \exp\left\{ - \frac{r \|\bm{v}\|^2}{   f_1\left[    \frac{\|\bm{u}\|^2}{  f_2( \theta )|_{[0,\pi]} }  \right]  }  \right\} \text{d}G(r) \text{d}\bm{v}\\
                       &  =  &  \pi^{k/2}  \frac{1}{ \{f_2( \theta )|_{[0,\pi]}\}^{l/2} }  \int_{[0,\infty)}   \exp\left\{ -  \frac{\|\bm{\omega}\|^2}{4r}  f_1\left[    \frac{\|\bm{u}\|^2}{  f_2( \theta )|_{[0,\pi]} }  \right] \right\} \text{d}\tilde{G}(r),
\end{eqnarray*}
where   the last equality follows from Fubini's Theorem and   $\text{d}G(r) = r^{k/2} \text{d}\tilde{G}(r)$.  In addition,  for any $\bm{\omega}$, the mapping $$g_{\bm{\omega}}(\cdot) := \int_{[0,\infty)}   \exp\left\{ - \frac{  \|\bm{\omega}\|^2}{4r}  f_1(\cdot) \right\} \text{d}\tilde{G}(r)$$ is completely monotone  \citep{feller1966introduction}. Therefore, 
 $$     \text{C}_{\bm{\omega}} (\theta,\bm{u})   =   \pi^{k/2}  \frac{1}{ \{f_2( \theta )|_{[0,\pi]}\}^{l/2} }  g_{\bm{\omega}}\left(   \frac{\|\bm{u}\|^2}{  f_2( \theta )|_{[0,\pi]} }    \right),$$
 and by Lemma \ref{prop_latent1},  we have that $\text{C}_{\bm{\omega}}\in \Upsilon_{d,l}$, for all $d\in\mathbb{N}$.

\footnotesize
\bibliographystyle{apalike}
\bibliography{mybib}

\end{document}